\documentclass{article}

\begin{document}
\newtheorem     {Def}              {Definition}    [section]
\newtheorem     {Lemma}    [Def]   {Lemma}
\newtheorem     {Theorem}  [Def]   {Theorem}
\newtheorem     {DefSatz}  [Def]   {Definition und Satz}
\newtheorem     {DefLemma} [Def]   {Definition und Lemma}
\newtheorem     {Bem}   [Def]   {Bemerkung}
\newtheorem     {Beh}   [Def]   {Behauptung}
\newtheorem     {Kor}   [Def]   {Korollar}
\newtheorem     {Bsp}   [Def]   {Beispiel}
\newtheorem     {Aufg}  [Def]   {Aufgabenstellung}
\newtheorem     {Cor}   [Def]   {Corollary}

\newenvironment {prqed}        {\begin{description}\item[Proof:]}{\newline\nopagebreak\hfill q.e.d. \end{description}\medskip}
     
\renewcommand   {\phi}  {\varphi}
\newcommand     {\R}    {{\ifmmode{\rm I}\mkern-4mu{\rm R}
    \else\leavevmode\hbox{I}\kern-.17em\hbox{R}\fi}}
\newcommand     {\C}    {\ifmmode{{\rm C}\mkern-15mu{\phantom{\rm t}\vrule}}
   \mkern10mu
   \else\leavevmode\hbox{C}\kern-.5em\hbox{I}\kern.3em\fi}
\newcommand     {\N}    {\ifmmode{\rm I}\mkern-3.5mu{\rm N}
    \else\leavevmode\hbox{I}\kern-.16em \hbox{N}\fi}
\newcommand     {\Z}    {\Bbb{Z}}
\newcommand     {\HR} {Hilbert\-raum}
\newcommand     {\Cz} {\C\,[ z\, ]}
\newcommand     {\Czzq} {\C\,[ z,\overline{z}\,]}
\newcommand     {\zq} {\overline{z}}
\newcommand     {\ggT} {\mbox{ggT}\,}
\renewcommand   {\H} {{\cal H}}  
\newcommand     {\eps} {\varepsilon}
\newcommand     {\Q} {{\cal Q}}
\newcommand     {\A} {{\cal A}}

\renewcommand   {\theenumi}   {(\roman{enumi})}
\renewcommand   {\theenumii}  {\alph{enumii})}
\renewcommand   {\theenumiii} {(\Roman{enumiii})}
\renewcommand   {\theenumiv}  {\Alph{enumiv})}
\renewcommand   {\labelenumi}   {(\roman{enumi})}
\renewcommand   {\labelenumii}  {\alph{enumii})}
\renewcommand   {\labelenumiii} {(\Roman{enumiii})}
\renewcommand   {\labelenumiv}  {\Alph{enumiv})}

\hyphenation{}

\begin{titlepage}
\begin{center}
{\Large
Orthonormal bases of polynomials in one complex variable 
}\bigskip\\
D.P.L. Castrigiano\\
W. Klopfer
\smallskip\\
Zentrum Mathematik der Technischen Universit\"at M\"unchen\\
D-80333 M\"unchen, Arcisstra\ss e 21, Germany

\smallskip

{\tt
castrig@mathematik.tu-muenchen.de\\
klopfer@mathematik.tu-muenchen.de}
\medskip\\
\small
\begin{quote}
Let a sequence $(P_n)$ of polynomials in one complex variable satisfy a recurrence relation with length growing slowlier than linearly. It is shown that $(P_n)
$ 
is an orthonormal basis in $L^2_{\mu}$ for some measure $\mu$ on $\C$, if and only if the recurrence is a $3-$term relation with special coefficients. The support of $\mu$ lies on a straight line. This result is achieved by the analysis of a formally normal irreducible Hessenberg operator with
only finitely many nonzero entries in every row. It generalizes the classical
Favard's Theorem and the Representation Theorem.

\end{quote}
\end{center}

\medskip

 \hspace*{22mm} recurrence relations, Hessenberg operator, spectral measure\\


Abbreviated title: orthogonal bases of polynomials
\end{titlepage}

                    \section{Restricted Recurrence}

Let $(P_n)_{n\geq 0}$ denote a sequence of polynomials in one complex variable
$z$ with $\deg P_n = n$ and $P_0 \equiv 1$. Such a sequence is uniquely  
determined by the {\bf recurrence relations} following by linear algebra
\begin{equation} zP_n (z) = \sum_{i=0}^{n+1} d_{in} P_i (z),\;\;\  d_{n+1,n}\not= 0
\label{eq1}
\end{equation}
for $n\ge 0$ with complex $d_{in}$. We call (\ref{eq1}) a {\bf restricted recurrence}  (rr) if there are $r_n\in\{0,...,n+1\}$ such that

\begin{equation} zP_n (z) = \sum_{i=r_n}^{n+1} d_{in} P_i (z)\quad {\mathrm and} \quad  r_n\rightarrow\infty.
\label{eq2}
\end{equation}
In other words the length of the recurrence $(\leq n+2-r_n)$ grows slowlier than linearly.\\
Examples of  rr  are the classical symmetric $3-$term recurrence where
$d_{nn}\in\R,\; d_{n+1,n} =\bar{d}_{n,n+1}$ for $n\geq 0$. By Favard's Theorem and the Representation Theorem (see e.g.[1; Chap.II,Th.6.4]         
or [2; 4.Cor.]) the corresponding sequences $(P_n)$ of
polynomials are exactly those which form an orthonormal basis (onb) in  
$L^2_\mu$ for some Borel  measure $\mu$  on $\C$  
with $\mu (\R
)=1$. The simplest rr is given
by   $d_{in}=\delta_{i,n+1}$ yielding the polynomials $P_n(z)=z^n$ which are  
orthonormal with respect to the normalized Lebesgue measure on the unit  
circle. This measure is unique and in particular there is no measure $\mu$ on $
\C $ such that the $z^n$ are an onb
in $L^2_\mu$. It can also be shown that for instance the polynomials $1,
2z, z^2, 2z^3,...$ satisfying an rr with $d_{in}=2^{(-1)^{n+1}}\delta_{i,n+1}$
do not admit any orthonormalizing measure at all [3; Beispiel 1.1.10]. Of
course, all recurrences of finite length are rr.\\
The main result we will obtain in (\ref{cor3.2}) is that an rr for which the polynomials $P_n$ yield an onb in $L^2_\mu$ for some measure $\mu$ on $\C$ is a $3-$term recurrence with special coefficients. The orthonormalizing measures are
concentrated on some fixed straight line. By an affine map of the complex
plane the polynomials can be transformed into orthogonal polynomials on
the real line.\\
The results are achieved by an operator theoretical approach. The main  
idea is to interpret (\ref{eq1}) as determining an irreducible Hessenberg operator in a Hilbert space, where the $P_n$ form an onb.\\                                                                                                               \section{Hessenberg Operator on \C[z]}  

The linear space $\C[z]$ of polynomials in one complex variable is equipped
with the scalar product $\langle \; ,\: \rangle  
$ by which the vector space basis $(P_n)$ is orthonormal. Let the Hilbert space ${\cal H}$ be a completion of $(\C[z],\langle\; ,\: \rangle )$. Then $\C [z]$ is the invariant dense domain of the operator $D$ in ${\cal H}$ given by
\begin{equation} (Dp)(z):= zp(z) \quad {\mathrm{for\; all\; polynomials}} \; p.
\label{eq3}
\end{equation}
With respect to the onb $(P_n)$, $D$ is a general {\bf irreducible Hessenberg
operator}; in particular no lower diagonal element $\langle
P_{n+1},DP_n \rangle =d_{n+1,n}$
vanishes. $D$ is a {\bf symmetric Jacobi operator} iff $(P_n)$ satisfies a
symmetric $3-$term recurrence. \\  
We begin with a rudimentary use of rigged
Hilbert space idea.\\

\begin{Lemma}                                         
  For $\lambda$ in $\C  
$ let $L$ be a linear form on $\C
[z]$ satisfying
$L(Dp)=\lambda                     L(p) \: \forall p\in \C
[z]$. Then there is a $c\!\in\! \C
$ such that
$L(p)=cp(\lambda)$. Hence $L=c\sum_n P_n(\lambda)\langle P_n,\cdot  
 \rangle$, since
 $p(z)=\sum_nP_n(z) \langle P_n,p \rangle $ holds  for all $z\in
\C $ and $p\in \C [z]$.
\label{le2.1}
\end{Lemma}
\begin{prqed} Since $\lambda L(P_n)=L(DP_n) \stackrel{(\ref{eq1})}{=}
L(\sum_{i=0}^{n+1}d_{in}P_i)=
\sum_{i=0}^{n+1}d_{in}L(P_i)$ holds, $(L(P_n))_{n\geq 0}$ satisfies
  the same recurrence relations (\ref{eq1})
 as $(P_n(\lambda))_{n\geq 0}$. Therefore $L(P_n)=cP_n(\lambda)$ holds
with $c:=L(P_0)$. By linearity it follows $L(p)=cp(\lambda)$. Finally note  
that the expansion $p=\sum_n \langle P_n,p \rangle P_n$ holds even pointwise, since the sum
is finite.
\end{prqed}
In the following analysis the relevant property of $D$ will be the formal
normality, i.e. ${\mathrm dom}(D^*)\supset\C[z]$ and $\| D^*p\| =
\| Dp\|\:\forall p\in\C[z]$.
Obviously $D$
is formally normal if it admits a normal extension in ${\cal H
}$. (We do not know
whether the converse is also true, cf. [4] or [3; Satz 4.3.14].)\\

\begin{Theorem}
 There is a bijective correspondence between the measures $\mu$ on
$\C$ for which $(P_n)_{n\geq 0}$ is an onb in $L^2_{\mu}$ and the normal
extensions $N$ of $D$ in ${\cal H}$. More precisely the following holds.\\
(i) If $(P_n)_{n\geq 0}$ is an onb in $L^2_{\mu}$ then ${\cal H}=L^2_{\mu}$
is a completion of $(\C[z],\langle\ ,\ \rangle )$ and the multiplication operator $M_z$ by
$z$ is a normal extension of $D$ in $L^2_{\mu}$.\\
(ii) Let $N$ be a normal extension of $D$ in ${\cal H}$. Denote by $E$ the
(projection valued) spectral measure of $N$ and let $\mu$ be the  
measure on $\C$ given by $\mu (\Delta):= \langle P_0,E(\Delta )P_0 \rangle $. Then the
support of ${\mu}$ equals the spectrum of $N$ and $(P_n)_{n\geq 0}$ is  
an onb in $L^2_{\mu}$. By the latter ${\cal H}$ and $L^2_{\mu}$ are identified
and $N$ equals $M_z$ (see (i)).
\label{th2.2}
\end{Theorem}
\begin{prqed} (i) is obvious. As to (ii) show first that $\langle E(\Delta)P_0:\;   \Delta\;{\mathrm  Borel \;  set} \rangle $ is dense in ${\cal H}$. Indeed, since the monomials $z^n=D^nP_0=N^n P_0$                     
 span $\C[z]$, they are total in ${\cal H}$. Furthermore $\infty > \|
 N^nP_0 \| ^2\; = \int | z^n |  ^2 d{\mu}$, whence $z^n \in L^2_{\mu}$. Therefore, for $\epsilon >o$, there is an elementary function $u=
\sum_i \alpha _i 1_{\Delta _i}$ such that $\epsilon \geq \int | u-z^n |
^2d\mu\;= \| (\sum_i \alpha _i E(\Delta _i)-N^n)P_0 \| ^2$.\

Therefore and since $\langle E(\Delta) P_0,E(\Delta ')P_0 \rangle = \langle P_0,E(\Delta \cap \Delta ') P_0 \rangle =\mu (\Delta \cap \Delta ')= \int \bar{1}_{\Delta }1_{\Delta '}d\mu$, the
assignment $E( \Delta) P_0\mapsto 1_{\Delta }$ determines a Hilbert space
isomorphism $\beta$ from ${\cal H}$ onto $L^2_{\mu}$. From $\beta E(\Delta)  
\beta ^{-1} 1_{\Delta '} = M_{1_{\Delta }} 1_{\Delta '}$, where
 $M_{1_{\Delta}}$
are the canonical projections of the spectral measure of $M_z$, it follows
$\beta N  
\beta ^{-1}=M_z$. This implies for $z^n \in \C [z] \subset {\cal H}$ that
$\beta z^n = \beta N^n P_0 = M_z^n \beta P_0 = M_z^n 1 =z^n \in L^2_{\mu}$.\\
Now (ii) follows easily. (ii) shows that $N$ and $N'$ are different if and
only if the corresponding measures $\mu$ and $\mu'$ are different. This sets
up the asserted bijective correspondence.  
\end{prqed}

\noindent
Consequently, there is a unique measure $\mu$ on $\C$ for which $(P_n)$ is
an onb in $L^2_{\mu}$, if $D$ is essentially normal. We do not know whether the converse is also true in general. Note however [3; Satz 4.2.6], where $(P_n)$ is not complete.
      
      \section{Formally normal Hessenberg operator in case of rr}

Because of $d_{k,k-1} \neq 0$ for $k>0$ no row of an irreducible Hessenberg operator vanishes, except for the $0$--th one (as e.g. in the case of the shift related to the polynomials $z^n$). But if $D$ is formally normal, then
$ \| D^*P_0 \| ^2 = \| DP_0 \| ^2 = \| d_{00}P_0 + d_{10}P_1 \| ^2 = |d_{00}|^2 + |d_{10}|^2 \geq |d_{10}|^2 > 0 $, so that there is $n \geq 0$ with $\langle P_0,DP_n \rangle = \langle D^*P_0,P_n \rangle \neq 0$. On the other hand, if $D$ represents an rr then for every $k>0$ there is an integer $s$ such that $d_{kn} = \langle P_k,DP_n \rangle = 0 $ for all $n>s$. This means equivalently that $\mbox{dom}(D^*) \supset \Cz$ and $D^* \Cz \subset \Cz $. In particular it follows by the polarization identity that, in case of an rr, $D$ is formally normal iff $D^*Dp = DD^*p \;\forall p \in \Cz$.
\begin{Theorem}
Let $D$ be formally normal representing an rr. Then
\[ D = b + aJ \]
where $b \in \C$, $a \in \C$ with $|a|=1$, and $J$ a symmetric Jacobi operator.\label{th3.1}
\end{Theorem}
\begin{prqed}
By the foregoing considerations for every $k \geq 0$ there is an integer $s_k \geq 0$ such that $d_{ks_k} \neq 0$ and $d_{kn} = 0$ for  $n > s_k$.
\smallskip\\
\setcounter{enumi}{1}
\theenumi
We are going to show that $s_k = k+1$ for all $k \geq 0$. For $\lambda \in \C$ denote by $F_{\lambda}$ the linear form $\sum\limits_n P_n(\lambda) \langle P_n, \cdot \rangle$ on $\Cz$, see (\ref{le2.1}). Because of $F_{\lambda}(D^*Dp) = F_{\lambda}(DD^*p) = \lambda F_{\lambda} (D^*p)\; \forall p \in \Cz$, $F_{\lambda} \circ D^*$ is a linear form on $\Cz$ satisfying the premise of (\ref{le2.1}). Therefore there is a $q(\lambda) \in \C$ such that $F_{\lambda} \circ D^* = q(\lambda)F_{\lambda}$. Again, by (\ref{le2.1}) it follows for $k \geq 0$ that\begin{eqnarray*}
q(\lambda) P_k(\lambda) & = & q(\lambda) F_{\lambda} (P_k) = F_{\lambda} (D^*P_k) = \sum\limits_n P_n(\lambda) \langle P_n,D^*P_k \rangle = \\
& = & \sum\limits_n P_n(\lambda) \overline{d}_{kn} = \sum\limits_{n'} P_{n'}(\lambda) \overline{d}_{kn'}
\end{eqnarray*}
with $\max \{ 0,k-1 \} \leq n' \leq s_k$. This implies for $k=0$ that $q$ is a polynomial in $\lambda$ of degree $s_0$. For general $k \geq 0$ it follows that $s_k = s_0 + k$ by comparing the degrees of the polynomials.\\
In order to prove $s_0 = 1$ observe that by the above computation the smallest index $m$ with $\langle P_m,qP_k \rangle \neq 0$ is $k-1$ if $k \geq 1$, since $d_{k,k-1} \neq 0$. Another way of calculating $m$, which follows now, yields $m = k - s_0^2$ if $k \geq s_0^2$. This shows $s_0 = 1$ as asserted.\\
For an integer $r \geq 0$ and $k \geq rs_0$ the smallest integer $m$ with $\langle P_m,\lambda^rP_k \rangle \neq 0$ is $k - rs_0$. This follows easily by iterating the recurrence $\lambda P_n = d_{n-s_0,n}P_{n-s_0} + \dots + d_{n+1,n}P_{n+1}$ with $d_{n-s_0,n} \neq 0$, which holds true if $n \geq s_0$, since $s_{n-s_0} = n$. Now, since $q$ is a polynomial of degree $s_0$, the assertion follows from this result for $r = s_0$.
\smallskip\\
\setcounter{enumi}{2}
\theenumi
Since $s_k = k+1$, $D$ represents a 3-term recurrence relation
\[DP_0 = g_0P_0 + h_0P_1 \: , \; DP_n = f_nP_{n-1} + g_nP_n + h_nP_{n+1}
\]
for $n>0$ and $h_n \neq 0$ for $n \geq 0$. By a straight forward computation the equations
\begin{itemize}
\item
$|f_1|^2 = |h_0|^2 \: , \; |h_{n-1}|^2 + |f_{n-1}|^2 = |f_n|^2 + |h_n|^2$
\item
$g_{n-1} \overline{h}_{n-1} + f_n \overline{g}_n = \overline{g}_{n-1}f_n + \overline{h}_{n-1} g_n$\item
$f_n \overline{h}_n = \overline{h}_{n-1} f_{n-1}$
\end{itemize}
follow for $n \geq 1$ from $DD^*P_n = D^*DP_n$ for $n \geq 0$. The last ones imply $f_n = \delta \overline{h}_{n-1}$ with $\delta := f_1 / \overline{h}_0$ for $n \geq 1$. Using these relations the equations of the second row yield $g_n - g_{n-1} = \delta ( \overline{g}_n - \overline{g}_{n-1} )$. Independently, the equations of the first row show $|f_n| = |h_{n-1}|$. Hence $| \delta | = 1$ follows.\\
Now choose $a \in \C$ satisfying $a^2 = \delta$ and put $b := g_0$. It is easy to show that $J := \overline{a} (D - b)$ is a symmetric Jacobi operator, as asserted.
\end{prqed}
Plainly, every operator $D$ of the form $b+aJ$ as in (\ref{th3.1}) represents an rr and is formally normal. More than that, it has normal extensions in ${\cal H}$, since every symmetric Jacobi operator has selfadjoint extensions (see e.g. [2; sec. 3] or [3] ), all of them explicitly given by the von Neumann Theory on symmetric operators. Conversely, $A := \overline{a} (N-b)$ is a selfadjoint extension of $J$ in ${\cal H}$. Indeed, $A$ is normal extending $J$. Hence $A^* \subset J^*$ and $\mbox{dom}(A^*) = \mbox{dom}(A)$. Since $J \subset J^*$, it follows $J \subset A^*$, whence $A \subset J^*$ and therefore $A = A^*$. In particular it follows that $D$ has a unique normal extension iff $J$ is essentially selfadjoint. Hence, by e.g. [2; 2.Theorem] or [3], uniqueness holds iff $\sum_{n \geq 0} |P_n(z)|^2 = \infty$ for some (and hence for every) $z \notin b + a \R$. Now (\ref{th2.2}) yields the following


\begin{Cor}
Let $(P_n)_{n \geq 0}$ be a sequence of complex polynomials satisfying $\mbox{\/\rm deg} P_n = n$ and $P_0 \equiv 1$. Then the statements (i) - (iv) are equivalent:
\begin{enumerate}
\item
$(P_n)$ satisfies an rr (\ref{eq2}) and is an onb in $L^2_{\mu}$ for some measure $\mu$ on $\C$.\item
$(P_n)$ satisfies a $3-$term recurrence of the kind
\[
zP_n(z) = a \overline{c}_{n-1}P_{n-1}(z) + (b+ad_n)P_n(z) + ac_nP_{n+1}(z)
\]
for $n \geq 0$ with $a,b,c_n \in \C$, $c_n \neq 0$, $|a|=1$, $d_n \in \R$, and $P_{-1} := 0$.
\item
The polynomials $\tilde{P}_n(z) := P_n(b+az)$, $n \geq 0$, for some $a,b \in \C$, $|a|=1$ satisfy a symmetric $3-$term recurrence\[
z \tilde{P}_n(z) = \overline{c}_{n-1} \tilde{P}_{n-1}(z) + d_n \tilde{P}_n(z) + c_n \tilde{P}_{n+1}(z) \: .
\]
(and hence form an onb in $L^2_{\nu}$ for some measure $\nu$ on the real line).
\item
$(P_n)$ is othonormal in $L^2_{\mu}$ for some measure $\mu$ on a line $b+a \R$ with $a,b \in \C$, $|a| = 1$.
\smallskip
\end{enumerate}Moreover, the values of $a,b$ in (ii) - (iv) are the same. Every orthonormalizing measure is concentrated on $b + a \R$; it is unique iff $\sum_{n \geq 0} |P_n(z)|^2 = \infty$ for some (and hence for every) $z \notin b + a \R$.
\label{cor3.2}
\end{Cor}

\begin{thebibliography}{1}
\bibitem{1}
T.S. Chihara, "An Inroduction to Orthogonal Polynomials", Gordon \& Breach, New York, 1978
\bibitem{2}
D.P.L. Castrigiano, Orthogonal Polynomials and Rigged Hilbert Space, J. Func. A. \underline{65}, (1986), 309--313
\bibitem{3}
W. Klopfer, Unbeschr\"ankte normale Operatoren im Hilbertraum und deren Anwendungen auf orthogonale Polynome, Thesis, Technical University Munich 1998
\bibitem{4}E.A. Coddington, Formally normal operators having no normal extension, Canad. J. Math. \underline{17}, (1965), 1030--1040
\end{thebibliography}
\end{document}